\documentclass[a4paper,twoside,final]{article}
\usepackage{amssymb}
\usepackage[latin1]{inputenc}
\usepackage[T1]{fontenc}
\usepackage{lpretex}
\usepackage[all]{xy}
\usepackage{amscd}
\usepackage{title}

\DocVersion[submitted]{1}
\DocVersion[minor correction and ref to Lazard]{2}

\candef s
\def\Num{\symb{\Cal N\!\mathit{um}}}
\def\Sets{\Cal S\!\mathit{ets}}
\def\Fr{\ensuremath{\Cal F\!\mathit{ree}}}
\def \Hom {\operatorname{Hom}}
\def \Aut {\operatorname{Aut}}
\def\sG{\Cal G}
\let\tens\otimes
\newcommand\spres[1]{\langle#1\rangle}
\newcommand\sprest[1]{{\langle#1\rangle}_t}
\begin{document}

\title[Integral minimal models]{On minimal models in integral homotopy theory}
\author[Torsten Ekedahl]{Torsten Ekedahl}
\address{Department of Mathematics\\
 Stockholm University\\
 S-106 91  Stockholm\\
Sweden}
\email{teke@matematik.su.se}

\maketitle

Inspired by some ideas of A.~Grothendieck (\cite{grothendieck83::en}) I will in
this article give an algebraic description of nilpotent homotopy types with
finitely generated homology (in each degree).

From a homotopy theoretic perspective the main result says that every such
nilpotent homotopy type may be represented by a simplicial set which is of the
form $\Z^n$ for some $n$ in each degree and for which the face and
degeneracy maps are \emph{numerical} maps, maps that are given by polynomials with
rational coefficients. Furthermore, the cohomology may be computed using
numerical cochains, any map between such models is homotopic to a numerical one
and homotopic numerical maps are numerically homotopic. The construction of a
model is rather straightforward; one first shows that the cohomology of an
Eilenberg Mac-Lane space can be computed using numerical cochains and then uses
induction over a principal Postnikov tower.

Localisation and completion fits very nicely into this framework. If $R$ is
either a subring of $\Q$ or is the ring of $p$-adic numbers for some prime $p$
and $K := R\Tensor\Q$ then a numerical function $\Z^m \to \Z^n$ that clearly
induces a map $K^m \to K^n$ takes $R^m$ to $R^n$. Hence a numerical simplicial
set gives rise to a numerical set obtained by replacing each $\Z^n$ that appears
in some degree by $R^n$. For a model this new space is the $R$-localisation when
$R$ is a subring of $\Q$ and the $p$-completion when $R=\Z_p$. In this way a
numerical model might be thought of as a universal localisation.

If $G$ is a finitely generated nilpotent torsion-free group, the theory may be
applied to $K(G,1)$ but stronger results are available. In fact $G$ has a
canonical structure of numerical group and for that structure the numerical
functions $G \to \Z$ are exactly the polynomial functions in Passi's sense (cf.,
\cite{passi68::dimen}). Furthermore, the cohomology of the group (with
$\Z$-coefficients) may be computed using numerical cochains (which are the same
as the Passi polynomial maps).

The theory can be reformulated in terms of cosimplicial rings. The cosimplicial
ring of numerical mappings of a model into $\Z$ has the property of being a free
\emph{numerical ring} in each degree, where a numerical ring is a ring together
with certain extra operations. The cosimplicial ring of all cochains is a
numerical ring and the numerical model can be obtained by the usual construction
of a free cosimplicial object homology equivalent to a given one. As a further
motivation for the relevance of numerical rings for homotopy we also note that
the cohomology of cosimplicial numerical rings admit an action of all cohomology
operations.

When passing to the rational localisation of a numerical model, the theory
should be compared to the theories of Quillen (\cite{Qu69}) and Sullivan
(\cite{sullivan77::infin}) of rational homotopy. In the case of Quillen's theory
the first step in his construction of a differential Lie algebra model is to
represent a nilpotent homotopy type by the simplicial classifying space of a
simplicial group $G$ that in each degree is the quotient of a finitely generated
free group by some element of the descending central series. That means that $G$
is a torsion free finitely generated nilpotent group and hence $K(G,1)$ is a
numerical space. As for Sullivan's approach the closest connection I have found
is by considering his spatial realisation of a differential graded algebra
model. We will define a natural quotient of that spatial realisation that has a
natural structure of $\Q$-numerical model.

To return to the starting point for the results of this section, Grothendieck's
idea was to represent any (locally finite) homotopy type by a simplicial set
that is $\Z^n$ and for which the face and degeneracy maps are polynomial
maps. While this seems reasonable for rational homotopy types as rational
homotopy theory is fairly linear (and is a consequence of the results of this
article) several people have pointed out problems with that idea (though I am
not aware of any proof that it is impossible). It is the suggestion of the
present article that polynomial maps should be replaced by numerical maps that
are the same as polynomial maps over the rationals. Of course, numerical maps
have appeared previously in homotopy theory in connection with polynomial
functors but I do not know if there is any relations with that theory. 
\begin{section}{Preliminaries}

It will be convenient to state a few preliminary results in some generality. If
$\Cal A$ is an additive category for which all idempotents have kernels, then
the category of homological complexes in $\Cal A$, concentrated in non-negative
degrees, is equivalent, by the Dold-Puppe constructions, to the category of
simplicial objects in $\Cal A$. Following \[May92], which will be our general
reference concerning simplicial results, we will use $N(-)$ for the functor
going from simplicial objects to complexes and $\Gamma(-)$ for the functor going
the other way. Similarly for cosimplicial objects and cohomological complexes
concentrated in non-negative degrees.

If $\Cal C$ now is a category having finite products (and in particular a final
object $*$), then if $X$ is a simplicial object in $\Cal C$ we say that $X$ is a
Kan complex in $\Cal C$ if the simplicial set $X^K$, where $(X^K)_n :=
\Hom_{\Cal C}(K,X_n)$, is a Kan complex for all $K\in {\cal C}$. If $X,Y\in
\s{\cal C}$, $\s{\Cal C}$ being the simplicial objects of $\s{\Cal C}$, then we
can define the function complex $Y^X$ as a simplicial set using the definition
in terms of $(p,q)$-shuffles as in \cite[6.7]{May92}. Then \cite[6.9]{May92}
goes through, so that if $Y$ is Kan, then so is $Y^X$. Furthermore, we say that
a sequence $F \to X \to Y$ of simplicial objects is a Kan fibration if, for all
$K\in {\cal C}$, $F^K \to X^K \to Y^K$ is a Kan fibration in the usual
sense. Finally, still following \cite[18.3]{May92}, we define for $F,B\in
\s{\cal C}$, $G$ a group object in $\s{\cal C}$, and a group action $G\times F
\to F$m a twisted cartesian product (\Definition{TCP}) to be an $E({\gta})\in
\s{\cal C}$ s.t.~$E({\gta})_n = F_n\times B_n$ and
\begin{eqnarray*}
\partial_i(f,b) =& (\partial_if,\partial_ib)\;      i>0,&(i)\\
\partial_0(f,b) =& (\tau(b)\cdot \partial_0f,\partial_0b), &(ii)\\
s_i(f,b) =& (s_if,s_ib),           &(iii)
\end{eqnarray*}
where \pil{\tau}{B_q}{G_{q-1}} is a morphism fulfilling the identities of
(\cite{May92}). If $F=G$ with the action being translation we will, still
following (\cite{May92}), speak of a principal twisted tensor product
(\Definition{PTCP}). It is then clear that if \pil T{\cal C} {{\cal C}'} is a
product preserving functor, then it takes TCP's to TCP's and PTCP's to
PTCP's. In particular, if $K\in {\cal C}$ then $(F^K,B^K,G^K,E(\tau)^K)$ is a
TCP and so, \cite[18.4]{May92}, $F \to E(\tau) \to B$ is a Kan fibration if $F$ is a
Kan complex and, in particular, $E(\tau)$ is Kan if $B$ is.

Suppose now that $A$ is a ring object in $\cal C$ s.\ t.\ multiplication induces
\begin{displaymath}
\Hom_{\cal C}(X\times Y,A)= 
\Hom_{\cal C}(X,A)\bigotimes_{\Z}\Hom_{\cal C}(Y,A)
\end{displaymath} 
for all $X,Y\in {\cal C}$ and that $\Hom_{\cal C}(X,A)$ is torsion free for all
$X\in {\cal C}$. If $F\in s{\cal C}$, $\Hom_{\cal C}(F,A)$ is a cosimplicial
abelian group and we put $H^*_A(X) := H^* (N(\Hom_{\cal C}(F,A)))$ and more
generally $H^*_A(X,M) := H^* (N(\Hom_{\cal C}(F,A))\bigotimes_{\Z}M)$, for $M$
an abelian group. Now, Szczarba's proof of the simplicial Brown's theorem \[Szc61]
uses only universal expressions and hence goes through in this context. If
$\tau$ is 1-trivial (i.e., $\tau_{|Bi}=*$, $i=0,1$) we therefore get a spectral
sequence
\begin{equation}\label{1.1}
E_2^{p,q} =  H_A^{*,p}(B,H_A^{*,q}(F))  \Rightarrow  H_A^{p + q}(E(\tau)).
\end{equation}
This is functorial for product preserving functors \pil T{\cal C}{{\cal C}'} and
this is so also if ${\cal C}'=\Sets$ and we make no particular assumption on
$T(A)$.

\end{section}

\begin{section}{Numerical spaces}

Having dealt with some preliminaries we can start getting down to business.
\begin{definition}
A \Definition{numerical function} from ${\Z}^m$ to ${\Z}^n$, $m,n\ge 0$, is a
function ${\Z}^m \to {\Z}^n$ which can be expressed by polynomials with rational
coefficients. If $M$ and $N$ are free abelian groups of rank $m$ resp.~$n$, then
a numerical function from $M$ to $N$ is a function $M \to N$ such that for one
(and hence any) choice of group isomorphisms ${\Z}^m \riso M$ and $N \riso
{\Z}^n$, the composite ${\Z}^m \riso M \to N \riso {\Z}^n$ is a numerical
function. The category $\Num$ has as objects finitely generated free abelian
groups and as morphisms the numerical functions.
\end{definition}
\begin{remark}
In the literature numerical functions sometimes appear under the
name of polynomial maps. I dislike this terminology as it could easily be
interpreted to mean maps given by polynomials with \emph{integer} coefficients,
indeed to me it seems the natural interpretation. Furthermore, it does not seem
unlikely that polynomial maps will have a role to play in homotopy theory and
hence deserve a name of their own.
\end{remark}
It is clear that $M,N \mapsto M\times N$ is a product in ${\Num}$ and that the
zero group is a final object so that ${\Num}$ has finite products. Furthermore,
addition makes all objects in ${\Num}$ abelian group objects and addition and
product makes $\Z$ a ring object. If $X\in s({\Num})$ then we will put
$H^*_{\Num}(X,{\Z}):=H_{\Z}^*(X)$. Define ${\Num}_i := \Hom_ {\Num}({\Z}^i,{\Z})
$ which thus is a ring, the ring of \Definition{numerical functions in $i$
variables}.
\begin{proposition}\label{2.2}
\part[i] ${\Num}_1 = \Dsum _{k\ge 0} {\Z}{x\choose k}$ where $x\choose k$
is the numerical function $n\mapsto {n\choose k}$.

\part[ii] ${\Num}_i = {\Num}_1\bigotimes {\Num}_1\bigotimes   \dots  \bigotimes {\Num}_1$ ($i$ times).

\part[iii] For $M,N\in {\Num}$,
\begin{displaymath}
\Hom_{\Num}(M\times N,{\Z})= 
\Hom_{\Num}(M,{\Z})\bigotimes_{\Z}\Hom_{\Num}(N,{\Z}).
\end{displaymath}

\part[iv] 
${\Num}$ is anti-equivalent to the full subcategory of the category of rings
whose objects are the rings isomorphic to ${\Num}_i$ for some $i$.
\begin{proof}
\DHrefpart{i} and \DHrefpart{ii} are well known. (Use $\Delta f(x) := f(x+1)-
f(x)$ and induction on the degree for \DHrefpart{i} and
${\gDe}_{last\;variable}$ for \DHrefpart{ii}.) Then \DHrefpart{iii} follows from
\DHrefpart{ii} as any object in ${\Num}$ is isomorphic to some ${\Z}_i$. As for
\DHrefpart{iv}, $M\mapsto \Hom_{\Num}(M,{\Z})$ gives a functor in one direction
and properly interpreted $R\mapsto \Hom_{Rings}(R,{\Z})$ will be a
quasi-inverse. Now, evaluation gives a natural morphism
\begin{displaymath}
\pil{{\phi}_M} M {\Hom_{Rings}(\Hom_{\Num}(M,{\Z}),{\Z})}
\end{displaymath}
in the category of sets. I claim that this map is a bijection. Indeed, we may
assume that $M={\Z}^i$ for some $i$ and using \DHrefpart{ii} that $i=1$. As
already $x={x\choose 1}$ separates points $\phi:=\phi_{\Z}$ is injective.  Since
${\Num}_1\bigotimes _{\Z}{\Q}= {\Q}[x]$, any ring homomorphism ${\Num}_1 \to
{\Z}$ is determined by what it does to $x$ and so $\phi$ is surjective. We now
want to show that any ring homomorphism $\Hom_{\Num}(N,{\Z}) \to
\Hom_{\Num}(M,{\Z})$ induces, through $\phi_M$ and $\phi_N$, a numerical map $M
\to N$. We reduce to $M={\Z}_i$ and $N={\Z}$ so we want to show that if
\pil{\rho}{{\Num}_1}{{\Num}_i} is a ring homomorphism then ${\Z}^i\owns n\mapsto
\rho(x)(n)\in {\Z}$ is a numerical function which is obvious by the definition
of ${\Num}_i$. Hence $R\mapsto \Hom_{Rings}(R,{\Z})$ maps the subcategory of
rings isomorphic to some ${\Num}_i$ into $\Num$ and by what we have just proved
it is a quasi-inverse to $M\mapsto \Hom_{\Num}(M,{\Z})$.
\end{proof}
\end{proposition}
The proposition immediately gives the following corollary:
\begin{corollary}  
The category $\s({\Num})$ is anti-equivalent to the full subcategory of the
category of cosimplicial rings consisting of those cosimplicial rings $R.$ for
which $R_n$ is isomorphic to some ${\Num}_i$ for all $n$.
\begin{proof}
\end{proof}
\end{corollary}
\begin{remark} 
It is this corollary that will allow us to interpret the results of 
this section in terms of cosimplicial rings and eventually an algebraic 
description of nilpotent homotopy types.
\end{remark}
Let us now note that if $M,N\in\Num$, then any group homomorphism $M \to N$
is a numerical function so that the category $\Fr$ of free abelian groups of
finite rank and homomorphisms embeds naturally in ${\Num}$. This is an additive
category where all idempotents have kernels so the equivalences $\Gamma$ and $N$
are defined. From the explicit description of $\Gamma$, \cite[p.\ 95]{May92}, it follows
that there is a unique function \pil T{{\N}^{\N}}{{\N}^{\N}} such that if $\dots
\to C_n \to C_{n-1} \to \dots $ is a complex in $\Fr$ and $r_{C_*}\in {\N}^{\N}$
is defined by $r_{C_*}(i)=\rk C_i$, then $T(r_{C_*})(i)= \rk \Gamma(C_*)_i$. If $M$
is a finitely generated abelian group we let $g(M)$ be the minimal number of
generators of $M$. If $X$ is a nilpotent space with $H_i(X,{\Z})$ finitely
generated for all $i\ge 0$, we put
\begin{eqnarray}\label{2.5}
g_n(X) &:=& \sum _{i=0}^\infty g(\pi_n^i(X)/\pi_n^{i+1}(X)),\\
h_n(X) &:=& \sum _{i=0}^\infty g({\rm torsion}(\pi_n^i(X)/\pi_n^{i+1}(X))),
\end{eqnarray}
where $\pi^0_n (X) := \pi_n(X)$ and $\pi_n^{i +1}(X) := \langle
\gamma(x)x^{-1}|x\in \pi_n^i (X);\;\gamma\in \pi_1(X)\rangle$.  We note also
that the rank of a f.\ g.\ free abelian group is invariant under numerical
isomorphisms so we may unambigously speak of its rank as an object in $\Num$.
\begin{definition}  
\part A \Definition{numerical space} is a connected numerical simplicial object
$X$ s.\ t.\ if $F$ is the forgetful functor ${\Num} \to \Sets$, then the natural
map $H^*_{\Num}(X,{\Z}) \to H^*(F(X),{\Z})$ is an isomorphism and the homology
groups $H_*(F(X),\Z)$ are finitely generated.

\part A \Definition{fibred numerical space} is a connected simplicial object $X$ 
which is an inverse limit
\begin{displaymath}
X= \ili(\dots \to X^{n+1} \to X^{n} \to \dots \to X^{0}.)
\end{displaymath}
where
\begin{enumerate}
\item $X^n_0$ is a point for all $n$.

\item For all $j$ there is an $N$ s.t.~$X_j^{n+1} \to X_j^{n}$ is an isomorphism
for all $n\ge N$.

\item For all $n$, $X^{n+1} \to X^{n}$ is a PTCP (in ${\Num}$) with
fiber some $\Gamma(M^{n})$, $M^{n}$ an \Fr-complex.

\item  $X^{0}=*$.
\end{enumerate}

\part A \Definition{special numerical space} is a fibred numerical simplicial
object $X$ s.\ t., using the notations of the previous definition, there is for
each $n$ integers $(h_n,i_n)$ so that $H^i(M^{n})=0$ for $i \ne h_n$ and
$H^{h_n}(M^n)= \pi^{i_n}_{h_n}(X)/\pi^{i_n+1}_{h_n}(X)$. We will call the
sequence $\dots \to X^{n+1} \to X^{n} \to \dots$ a \Definition{special tower}.

\part  A special numerical space  $X$ is \Definition{minimal} if    
\begin{displaymath}
\forall n: \rk M^n = \sum_iT(\tau_i)(n),
\end{displaymath}
where
\begin{displaymath}
\tau_i(j)=\left\{            
\begin{array}{r@{,\quad\mathrm{ if\ }}l}
 g_i(X)& j=i,\\
 h_i(X)& j=i+1,\\
      0& j\ne i, i+1.
\end{array}
\right.
\end{displaymath}
We will call the sequence $\dots \to X^{n+1} \to X^{n} \to \dots$ a
\Definition{minimal tower}.
\end{definition}
\begin{remark} 
For $X$ a special numerical space
\begin{displaymath}
\forall n: \rk M^n \ge \sum_iT(\tau_i)(n),
\end{displaymath}
which justifies the terminology minimal.
\end{remark}
\begin{lemma}\label{2.7}  
\part[i]  If   $F \to E \to B$ is a  $1$-trivial  TCP in  $\s({\Num})$ and if
$F$ and  $B$ are numerical spaces then so is  $E$.

\part[ii]\label{2.7:ii} A special numerical space $X$ is a numerical space and a
Kan object. In particular, $F(X)=\Hom_{\Num}(*,X)$ is Kan.

\part[iii]  If  $X$ is a special numerical space then
\begin{displaymath}
\forall n: \rk X_n\ge \sum _iT(\tau_i)(n).
\end{displaymath}                   
\begin{proof}
Indeed, \DHrefpart{i} follows immediately from the map between the spectral
sequences of \(1.1) for the TCP and the induced TCP in $\s(\Sets)$. As for
\DHrefpart{ii}, the trivial nature of the limit in the definition of a special
numerical space allows us to assume that $X$ equals some $X_n$ and then by
\DHrefpart{i} that $X=\Gamma(M)$, for some \Fr-complex $M$. As $H^*(M)$ is
concentrated in one degree, $M$ is homotopic to a bounded complex $M'$, indeed
to one concentrated in 2 degrees, and then using the naive truncations of $M'$
and \DHrefpart{i} again we may assume that $M=N[n]$ for some object $N$ in
\Fr. Now there is a numerical (indeed a \Fr-) PTCP $\Gamma(N[n-1]) \to X \to
\Gamma(N[n])$ with $X$ \Fr-contractible so by induction on $n$ and Zeeman's
comparison theorem (cf., \[Ze57]) we may assume that $n=1$. We have a $\Fr$-PTCP
\begin{displaymath}
\Gamma(N_1[1]) \to \Gamma((N_1\times N_2)[1]) \to \Gamma(N_2[1]),
\end{displaymath}
so we may assume that $N={\Z}$. Hence we are reduced to showing that
$H^*_{\Num}(K({\Z},1)) \to H^*(K({\Z},1))$ is an isomorphism. A numerical
$1$-cocycle (for K({\Z},1)) is just a numerical function \pil f{\Z}{\Z}
s.t.~$f(x+y)=f(x)+f(y)$ so $f(x)=ax+b$ for some $a,b\in {\Z}$ and a 1-coboundary
is of the form $f(x)=c$. As we have exactly the same description for set
theoretical 1-cochains and 1-coboundaries we get an isomorphism for $*=1$. As
$*=0$ is trivial it only remains to show that $H^i_{\Num}(K({\Z},1))=0$ for
$i\ge 2$ as this is true in the set case. Now $K({\Z},1)_n={\Z}^{n+1}$ and all
the face operators are projections or sums of two adjacent coordinates. We can
grade ${\Num}_i$ by
\begin{displaymath}
\deg {x_1\choose n_1}{x_2\choose n_2}\dots {x_k\choose n_k}=\sum_in_i
\end{displaymath}
and then $\Hom_{\Num}(K({\Z},1),{\Z})$ becomes a graded complex as
\begin{displaymath}
{x+y\choose n} = \sum _{i+j=n}{x\choose i}{y\choose j}.
\end{displaymath}
 Hence to show the required vanishing we can replace
\begin{displaymath}
\Num_i = \sum \Z{x_1\choose n_1}{x_2\choose n_2}\dots {x_k\choose n_k}
\end{displaymath}
by
\begin{displaymath}
\ovl {\Num_i} = \prod {x_1\choose n_1}{x_2\choose n_2}\dots {x_k\choose n_k}
\end{displaymath}
 as the cohomology can be computed degree by degree. This gives us a 
complex $T$, say. Now, $\ovl {\Num_i} =  \Hom({\N}^i,{\Z})$ where
\begin{displaymath}
{x_1\choose n_1}{x_2\choose n_2}\dots {x_k\choose n_k}\mapsto 
\symb{char\ fct\ of}(n_1,n_2,\dots ,n_k)
\end{displaymath}                           
and this equality respects maps induced by projections and sums of coordinates
(but it is not a ring isomorphism). Therefore, $T$ is additively isomorphic to
the standard cochain complex of $K({\N},1)$ and so
$H^*(T)=H^*(K({\N},1),{\Z})=Ext^*_{{\Z}[t]}({\Z},{\Z})$ and the latter group is
clearly concentrated in degree 0 and 1. That a special numerical space is Kan
follows by induction on the Postnikov tower and a trivial passage to the
limit. Finally, to prove \DHrefpart{iii} it suffices to prove that if $M$ is a
$\Fr$-complex with a single non-zero homology group $H_i(M)$, then $\rk M_i \ge
g(H_i(M))$ and $\rk M_{i+1} \ge g(tor(H_i(M)))$. This, however, is clear by the
principal divisor theorem.
\end{proof}
\end{lemma}
\begin{remark}
\part The next to last part of the proof of ii) looks somewhat mysterious and
may be clarified by noting that ${\Num}_i$ is the ring of invariant differential
operators on the formal $i$-dimensional torus. Its coproduct is therefore dual
by Cartier duality \cite[II,\S4]{Dem72}) to the product on the coordinate ring
of the formal $i$-dimensional torus. Similarly, $\Hom({\N}^i,{\Z})$ is the
Cartier dual of the $i$-dimensional formal additive group. As the formal
$i$-dimensional torus and the $i$-dimensional formal additive group are
isomorphic as formal schemes but not as formal group schemes, ${\Num}_i$ is
isomorphic to $\Hom({\N}^i,{\Z})$ as coalgebras but not as rings. However, in
defining the differentials of chains on $K(-,1)$ only the coproduct is used.

\part Had we worked with polynomial instead of numerical functions everything
would have worked up to the statement $H^i_{pol}(K({\Z},1))=0$ for $i>1$. This
statement is false however. As a matter of well known fact, the polynomial 2-cocycle
$((x+y)^p-x^p-y^p)/p$, $p$ prime, is not the coboundary of a polynomial
1-cochain. It is the boundary of the numerical 1-cochain $(x^p-x)/p$ as the
lemma predicts.
\end{remark}
\begin{lemma}\label{2.8}
Let $Y$ be a special numerical space and $X$ a numerical space. Then $F$, the
forgetful functor, induces a bijection
\begin{displaymath}
\Hom(X,Y)/(\mathrm{numerical\ homotopy}) \to 
\Hom(F(X),F(Y))/(\mathrm{homotopy})
\end{displaymath}
where  $\Hom(-,-)$ means based maps.
\begin{proof}
Indeed, it will be easier to prove a stronger statement. Let $Z,V\in {\Num}$ be
pointed objects and $V$ Kan and put $[Z,V]_n:=\pi_n((V^Z)_b,*)$, where $(-)_b$
denotes based maps. As above, if $V_1 \to V_2 \to V_3$ is a TCP of Kan objects
then $(V^{Z_1})_b \to (V^{Z_2})_b \to (V^{Z_3})_b$ is a fibration and we get a
long exact sequence of homotopy. If $V_1 \to V_2 \to V_3$ is a PTCP then we get
as usual the extra precision that the fibers of $[Z,V_2]_0 \to [Z,V_3]_0$ are
the orbits under an action of $[Z,V_1]$ and the sequence extends to $[Z,V_2]_0
\to [Z,V_3]_0 \to [Z,\ovl WV_1]_0$ (cf.~\cite[p.~87]{May92}). We now consider a
sequence of PTCP's $Y^n \to Y^{n-1}$ as in the definition of a special numerical
space. We want to prove by induction on $n$ that $[X,Y^n]_m \to [F(X),F(Y^n)]_m$
is a bijection for all $m$. The case $n=0$ certainly causes no problem and in
general we have the PTCP $Y^n \to Y^{n-1}$ with fiber some $\Gamma(M^n)$. The
extra precision given to the long exact sequence is exactly what is needed to
make the 5-lemma work and we reduce hence to showing that $[X,\Gamma(M^n)]_m \to
[F(X),F(\Gamma(M^n))]_m$ is a bijection. Now, $M^n$ is homotopic to a bounded
complex, so we may assume that $M^n$ is bounded. By the same dévissage as before
we reduce to $M^n={\Z}[0]$ and then this bijection is true by the definition of
numerical space. Putting $m=0$ we then get the lemma for $Y$ replaced by
$Y^n$. To pass from $Y^n$ to $Y$ we use the Milnor exact sequence
\begin{displaymath} 
*\to \ili^1[X,Y^n]_1 \to [X,Y] \to \ili[X,Y^n] \to *,
\end{displaymath} 
the similar sequence for $F(X)$ etc and the $5$-lemma.
\end{proof}
\end{lemma}
We have now come to the main result of this section.
\begin{theorem}\label{set => num}
Let $X$ be a simplicial set which is nilpotent (which to us will include being
connected) of finite type (i.\ e., finitely generated homology in each degree).

\part  There is a minimal special numerical space  $Y$ and a homotopy 
equivalence  $X \to F(Y)$.

\part  If  $Y'$ is a numerical space and  $X \to F(Y')$ a homotopy equivalence, 
there is a unique, up to numerical homotopy, numerical morphism  
$Y' \to Y$, where  $Y$ is as in  i), making the following diagram commute
up to homotopy
\begin{displaymath}
\diagram{                              
&& F(Y)\cr
&\nearrow&\uparrow\cr
  X& \to  &F(Y').\cr}
\end{displaymath}

\part If  $Y'$ is also special then  $Y' \to Y$ is a numerical homotopy 
equivalence.

\begin{remark} 
It is not true that a homology equivalence between minimal 
special numerical spaces is necessarily an isomorphism as is shown by 
the following example:
\begin{displaymath}
\Gamma({\Z}\mapright2{\Z})\mapright3\Gamma({\Z}\mapright2{\Z}).
\end{displaymath}
Indeed, multiplication by $3$ induces an isomorphism on the homology, $\Z/2$,
but is clearly not an isomorphism.
\end{remark}
\begin{proof}
Let $\dots \to X^n \to X^{n-1} \to \dots $ be a minimal principal Postnikov
system, i.e., the $X^n \to X^{n-1}$ are principal fibrations with fiber some
$K(\pi_m^i (X)/\pi_m^{i+1}(X),m)$ using the notations of \(2.5) in order of
increasing $i$ and $m$. We will step by step replace, up to homotopy, $X^n$ by
some $F(Y^n)$ and $X^n \to X^{n-1}$ by $F(-)$ of a PTCP $Y^n \to Y^{n-1}$. Assume
that we have done this up to $n-1$. We then have a cartesian diagram, where
$K(M,m)$ is the fiber of $X^n \to X^{n-1}$,
\begin{displaymath}
\diagram{
&&&&K(M,m)\cr
&&&&\downarrow\cr
&&X^n & \to  &T\cr
&&\downarrow &&\downarrow \cr
F(Y^{n-1})&\riso & X^{n-1}& \to &K(M,m-1).\cr
}
\end{displaymath}
Here the right hand fibration is the standard one with $T\sim*$. Choose a
resolution $F_* \to M$ by free f.g.~abelian groups s.t.~$\rk F_0=g(M)$, $\rk
F_1=g(\symb{tor}M)$ and $F_i=0$ for $i>1$. There is then a numerical PTCP
$\Gamma(F.[m]) \to I \to \Gamma(F.[m+1])$ s.t.~$F(-)$ applied to it is homotopic
to $K(M,m) \to T \to K(M,m+1)$. Hence by lemma \ref{2.7} there is a morphism
\pil{\rho} {Y^{n-1}}{\Gamma(F.[m+1])} such that $F(\rho)$ is homotopic to
$Y^{n-1} \to X^{n-1} \to K(M,m+1)$.  Let $Y^n \to Y^{n-1}$ be the PTCP induced
by $\rho$ from $\Gamma(F.[m]) \to I \to \Gamma(F.[m+1])$. Then $F(Y^n) \to
F(Y^{n-1})$ is homotopic to ${X^n \to X^{n-1}}$ and we put $Y := \ili(\dots \to
Y^n \to Y^{n-1} \to .\dots )$. Then there is a morphism $X \to F(Y)$ which by
construction is a homology equivalence and so a homotopy equivalence as $X$ and
$F(Y)$ are nilpotent. If $X \to F(Y')$ is a homology equivalence, where $Y'$ is
a numerical space, then by obstruction theory applied to the $F(Y^n) \to
F(Y^{n-1})$ there is a map $F(Y') \to F(Y)$ s.\ t.\
\begin{displaymath}
\diagram{  
&&F(Y')\cr
&\nearrow&\downarrow \cr
 X &  \to &F(Y)}
\end{displaymath}
commutes up to homotopy. By lemma \ref{2.8} there is a morphism $Y' \to Y$ 
s.t.~$F(Y') \to F(Y)$ is homotopic to the given $F(Y') \to F(Y)$. In case $Y'$ also is 
special another application of lemma \ref{2.8} shows that $Y' \to Y$ is a numerical 
homotopy equivalence.  
\end{proof}
\end{theorem}
As we will see in an example in the next section, minimal models are not unique
up to isomorphism. In the simple case we can however say that a minimal tower
must be preserved.
\begin{proposition}\label{Postnikov factoring}
Let $X$ and $Y$ be simple minimal numerical spaces and $X^{\cdot}$ and $Y^\cdot$ their
minimal towers. Any map \pil{f}{X}{Y} induces a map \pil{f^\cdot}{X^\cdot}{Y^\cdot}.
\begin{proof}
The proof is easily reduced to the following statement. If $X' \to X$ is a PTCP
for a simplicial group $\Gamma(M_.)$ with $M_i=0$ if $i \le n$ and $Y' \to Y$ is
a PTCP for a simplicial group $\Gamma(N_.)$ with $N^i=0$ for $i > n+1$ and for which
$N_{n+1} \to N_n$ is injective then for any commutative diagram
\begin{displaymath}
\begin{CD}
X' @>f'>> Y'\\
@VgVV   @VVV\\
X @>>> Y
\end{CD}
\end{displaymath}
there is a factorisation $X \to Y'$ of the diagram. This again amounts to saying
that $f'$ is constant on a fibre of $g$. In proving this we immediately reduce
to the case when $X=Y=\Delta^m$, and then need only show that the
map is constant on the fibre over $e_m$, the unique non-degenerate $n$-simplex
of $\Delta^n$. As PTCP's over a simplex are trivial we may assume that $X' \to
X$ and $Y' \to Y$ are trivial PTCP's. Hence we reduce to showing that any map
$\Gamma(M_.)\times\Delta^m \to \Gamma(N_.)$ is constant on
$\Gamma(M_.)_m\times\{e_m\}$. In general a map $\Gamma(M_.)\times\Delta^m \to
\Gamma(N_.)$ is the same thing as an additive map $\Z[\Gamma(M_.)\times\Delta^m]
\to \Gamma(N_.)$ and what we want to show is that $[(k,e_m)]-[(0,e_m)]$ is
mapped to zero for $k \in \Gamma(M_.)_m$.  If $m \le n$ then this is
obvious. When $m \ge n$ want to show that the kernel $K$ of the map
$\Z[\Gamma(M_.)\times\Delta^m] \to \Z[\Delta^m]$ induced by projection maps to
zero in $\Gamma(N_.)$. To prove this is equivalent to showing that $N(K) \to
N(\Gamma(N_.))=N_.$ is zero. Now, $N(K)$ is a complex of free abelian groups,
$H_i(M_.)$ is zero when $i > n$ and $K_i$ is zero when $i \le n$ so that any map
$N(K) \to N_.$ is null-homotopic. However, $M_i$ is zero when $i > n+1$ so any
homotopy is zero.
\end{proof}
\end{proposition}
\end{section}
\begin{section}{Localisation and completion}

We will now extend the theory presented so far to some other base rings than the
(sometimes only implicitly mentioned) ring of rational integers. Our choice of
rings is dictated on the one hand by which rings that are being used for
defining localisation and completion in homotopy theory, on the other hand by
which rings for which a straightforward generalisation of the notion of
numerical function admits a description similar to the one given for the
integers. Somewhat surprisingly these two requirements seem to give the same
answer.
\begin{definition-lemma}
Let $R$ the rational numbers or the ring $\Z_p$ of
$p$-adic integers. An ($R$-)numerical function $F \to G$ between finitely
generated free $R$-modules is a function that can be defined by polynomials with 
coefficients in $R\Tensor\Q$.

\part The $R$-algebra, $\Num_n(R)$, of numerical functions $R^n \to R$ is free
as $R$-module on the functions $(r_1,r_2,\dots,r_n) \mapsto {r_1\choose
m_1}\dots{r_n\choose m_n}$.

\part The evaluation map $R^n \to \Hom_{R-algebras}(\Num_n(R),R)$ is a bijection.
\begin{proof}
Let us first prove that a polynomial with rational coefficients mapping $\Z^n$
to $\Z$ will map $R^n$ to $R$. If $R$ is a subring of $\Q$ then it is an
intersection of the localisations $\Z_{(p)}$ that contains it so in that case
one is reduced to $R = \Z_{(p)}$ and then to $R = \Z_p$ as $\Z_{(p)} = \Q \cap
\Z_p$. Our polynomial defines a continuous function $\Z_p^n \to \Q_p$, $\Z_p
\subset \Q_p$ is a closed subset and $\Z^n \subset \Z_p^n$ is a dense subset. As
$\Z^n$ is mapped into $\Z \subset \Z_p$ it follows that $\Z_p^n$ is mapped into
$\Z_p$. Conversely, if we have a polynomial with $R\Tensor\Q$-coefficients
mapping $R^n$ into $R$ it maps in particular $\Z^n$ into $R$ and a slight
modification of the argument in the proof of proposition \ref{2.2} show that it
is an $R$-linear combination of products of binomial polynomials.

The proof of the second part is entirely analogous to the same statement for $R
= \Z$ given in the proof of proposition \ref{2.2}.
\end{proof}
\end{definition-lemma}
\begin{remark}
One may wonder whether torsion free rings $R$ other than the ones mentioned in
the lemma have the property that elements of $\Num_n$ maps $R^n$ into $R$. Of
course, $\Num_n$ itself or any ring containing $\Q$ is such an example but one
can show that if $R$ is a finitely generated ring and if $P(r) \in R$ for all $r
\in R$ and all numerical polynomials $P \in \Num_1$ then $R \subset \Q$. Indeed,
we have that $(x^p-x)/p \in \Num_1$ for all primes $p$ and hence $r^p-r \in
pR$. From this one concludes that $r$ is algebraic for all $r \in R$ and hence
that $R$ is a subring of a number field $K$. The condition that $r^p-r \in pR$
for all primes $p$ and $r \in R$ then implies that almost all primes in $K$ are
of degree $1$ which implies that $K=\Q$ (details of this type of argument can be
found in \cite{Ka70}).
\end{remark}
We will say that a ring $R$ is a \Definition{coefficient ring} if it is a
subring of $\Q$ or equals the ring $\Z_p$ of $p$-adic integers for some prime
$p$. By some abuse of language we will say that a nilpotent simplicial set is
\Definition{$R$-local} if it is $R$-local in the usual sense if $R \subseteq \Q$
and is $p$-complete if $R = \Z_p$. Similarly we will speak of the
$R$-localisation of a nilpotent simplicial set.

From the lemma one can continue almost verbatim and introduce, when $R$ is a
coefficient ring, (special) $R$-numerical spaces (the condition being that
cohomology with $R$-coefficients can be computed using $R$-numerical
cochains). One also, though we shall not use it, gets that one may represent any
$R$-local nilpotent finite type space by a special $R$-numerical space. We will however
note that if $F \to G$ is a numerical map between finitely generated free
abelian groups and $R$ is a subring of $\Q$ or $\Z_p$ then we get an induced
$R$-numerical map $F \Tensor R \to G \Tensor R$. This gives a functor, also denoted
by $- \Tensor R$ from simplicial numerical objects to simplicial $R$-numerical
objects. We also get a map of simplicial sets $F(X) \to F(X \Tensor R)$. For
special numerical spaces this is a localisation map:
\begin{theorem}\label{localisation}
Let $R$ be a coefficient ring.

\part A special $R$-numerical space is $R$-local.

\part If $X$ is a special numerical space then $F(X) \to F(X\Tensor R)$ is an
$R$-localisation map.
\begin{proof}
The first part uses the fact that locality is stable under
fibrations, that $K(M,n)$ is local if $M$ is a finitely generated
free $R$-complex which is clear as its homotopy groups are and that by
speciality one may reduce to such $K(M,n)$'s. For the second part we again
reduce to $K(M,n)$'s for $M$ a finitely generated free $\Z$-complex.
\end{proof}
\end{theorem}
\begin{remark}
I do not know if $X \Tensor R$ is a $R$-numerical space if $X$ is a numerical
space nor if it always is local.
\end{remark}
\begin{proposition}
Let $X$ and $Y$ be minimal $R$-numerical spaces where $R$ is $\Q$, $\Z_{(p)}$ or
$\Z_p$. Then any homotopy equivalence \pil{f}{X}{Y} is homotopic to a map that
is the inverse limit of a map of inverse systems $X^{\cdot} \to Y^{\cdot}$ (where $X^.$
and $Y^.$ are sequences as required in the definition of minimality) such that
each $X^n \to Y^n$ is an isomorphism at each point. In particular, a homotopy
equivalence between minimal numerical spaces is homotopic to an isomorphism and
even more particularly minimal models of the same space are isomorphic.
\begin{proof}
Let $p$ be the characteristic of $R$ modulo its maximal proper ideal.  For
evident reasons we will have to carefully distinguish between equality versus
homotopy of maps and we will start off with some observations. They will apply
equally well to simplicial sets as to numerical spaces but for simplicity we
will speak only of simplicial sets;N in any case the numerical case may be
deduced from the set-theoretic one using \(set => num). To begin with, if $G$ is
a simplicial abelian group with a single homotopy group $M$ in degree $n \ge 0$
then isomorphism classes of PTCP's with structure group $G$ over a simplicial
set $X_.$ correspond to elements of $H^{n+1}(X,M)$ (cf., \cite[class of
PTCP's]{May92}). It follows from that proof together with the use of the mapping
cone construction that if $X \to Y$ is a map of simplicial sets then the
relative cohomology $H^{n+1}(Y,X,M)$ correspond to equivalence classes of PTCP's
over $Y$ together with a trivialisation of its pullback to $X$ where two of them
are equivalent if they are isomorphic over $Y$ by an isomorphism whose pullback
to $X$ is homotopic to one that preserves the given trivialisations. Let us also
note that as we are dealing with principal fibrations, giving a trivialisation
is the same thing as giving a section.

The way the Postnikov tower $\{X^{\cdot}\}$ fits into this description is that
$X^{n+1} \to X^{n}$ is universal for PTCP's in degree $h_n$ over $X^{n}$ that
are provided with a trivialisation over $X$. From this we can construct the maps 
by induction over $n$. We therefore may assume we have the following diagram
that is assumed to commute up to homotopy:
\begin{displaymath}
\begin{CD}
X @>f>> Y\\
@VVV @VVV \\
X^{n+1} && Y^{n+1}\\
@VVV @VVV \\
X^{n} @>f^{n}>> Y^{n}\\
\end{CD}
\end{displaymath}
and $f^n$ is an isomorphism. Let $M_n$ resp.\ $N_n$ be the complexes for which
$X^{n+1} \to X^n$ resp.\ $Y^{n+1} \to Y^n$ are $\Gamma(M_n)$- resp.\
$\Gamma(N_n)$-PTCP's. As $H^n(M_\cdot)$ resp.\ $H^n(N_\cdot)$ are $\pi^{i_n}_{h_n}(X)$
resp.\ $\pi^{i_n}_{h_n}(Y)$, $f$ induces an isomorphism between them. We lift
this isomorphism to a map of complexes $M_. \to N_.$. Now, as $N_.$ is minimal,
the image of $N_{n+1}$ in $N_n$ is contained in $pN_n$ and hence by Nakayama's
lemma the map $N_n \to N_n$ is a surjection. As $M_.$ also is minimal, the rank
of $M_n$ is the same as that of $N_n$ and so the map $M_n \to N_n$ is an
isomorphism. This implies that $M_. \to N_.$ is an isomorphism.

Now, the pullback of $Y^{n+1} \to Y^n$ along the composite $X \to Y \to Y^n$ has
a section and hence a trivialisation. As the diagram is homotopy commutative we
get a trivialisation of the pullback of $Y^{n+1} \to Y^n$ along the composite $X
\to X^n \to Y^n$. This in turn, by the universality of $X^{n+1} \to X^n$, gives
a mapping from $X^{n+1} \to X^n$ to the pullback of $Y^{n+1} \to Y^n$ along $X^n
\to Y^n$ covering the isomorphism $\Gamma(M_.) \to \Gamma(N_.)$. This again is
nothing but a map \pil{f^{n+1}}{X^{n+1}}{Y^{n+1}} of PTCP's covering
$\Gamma(M_.) \to \Gamma(N_.)$. As the latter as well as the base map, $f^n$, are
isomorphisms so is $f^{n+1}$. By construction it gives rise to a homotopy
commutative diagram
\begin{displaymath}
\begin{CD}
X @>f>> Y\\
@VVV @VVV \\
X^{n+1} @>f^{n+1}>> Y^{n+1}
\end{CD}
\end{displaymath}
and thus finishes the induction step.
\end{proof}
\end{proposition}
\begin{remark}
Uniqueness of minimal models is not true over the integers: Fix an integer $n >
1$ and consider a class $\alpha$ of order $n$ in $H^4(K(\Z/n,1),\Z)=\Z/n$. This can be used
as $k$-invariant for a fibration over $\Gamma(\Z[2] \mapright{n} \Z[1])$ with
fibre $K(\Z,3)$. Now, multiplication by any invertible residue $\beta$ modulo $n$ on
$H^2(\Z/n,\Z)=\Z/n$ can be induced by a homotopy equivalence of the base. Taking
into account also the action of multiplication by $-1$ on $K(\Z,4)$ we see that
two $k$-invariants $\alpha$ and $\alpha '$ give homotopic total minimal models
if (and only if) $\alpha '=\pm \beta^2\alpha$.

On the other hand it follows from \(Postnikov factoring) that any isomorphism
between two such models induces an isomorphism over $\Gamma(\Z[2] \mapright{n}
\Z[1])$. Let us first consider the induced isomorphism on the base
$B:=\Gamma(\Z[2] \mapright{n} \Z[1])$. As $B$ consists of a point in degree $0$
any map from $B$ to itself preserves the base point $0$. Then in degree we have
a numerical isomorphism from $\Z$ to $\Z$ taking $0$ to $0$. This in turn is a
polynomial isomorphism from $\Q$ to $\Q$ and as such is well known to have the
form $x \mapsto ax+b$ and as $0$ is preserved $b=0$, as $\Z$ is preserved $a\in
\Z$ and as its inverse has the same properties $a =\pm1$. Now, it is easy to see
that a map $B \to B$ is determined by what it does in degree $1$ so any
automorphism of $B$ is given by multiplication by $\pm1$. Multiplication by $-1$
acts trivially on the $k$-invariants in question so we may assume that the
induced map on $B$ is the identity. For the rest of the argument we will ignore
the numerical structure. Any $K(\Z,3)$-PTCP over $B$ is classified as fibration
over $B$ by a torsor over the simplicial set of automorphisms of the simplicial
set $K(\Z,3)$. On the one hand we have the translations which is isomorphic as
simplicial set $K(\Z,3)$. Using them we may concentrate on based isomorphisms.
If we more generally consider the simplicial set of based endomorphisms of
$K(\Z,3)$ then the argument of \(Postnikov factoring) shows that they are
determined by the action on the third homotopy group so that the simplicial
group of based automorphisms is equal to the constant simplicial group
$\{\pm1\}$. Hence the simplicial group of automorphisms of $K(\Z,3)$ is the
split extension of $K(\Z,3)$ by the constant group $\{\pm1\}$ acting by
multiplication. From this it follows that two $K(\Z,3)$-PTCP's that are
isomorphic as fibrations either are isomorphic as PTCP's or one is isomorphic to
the transformation by multiplication of $-1$ on the other. In the notation above
that means the relation $\alpha ' = \pm1 \alpha$. Hence, there are in general
minimal models that are homotopic but not isomorphic.
\begin{remark}
Note that the fibrations we get are exactly the first non-trivial step in the
Postnikov tower of the three-dimensional lens spaces. I have no idea whether the
fact that the isomorphism classes of minimal models coincides with the
homeomorphism classes of these lens spaces has any significance.
\end{remark}
\end{remark}
\end{section}
\begin{section}{Cosimplicial ring interpretation}

We now want to interpret what we have proved in terms of cosimplicial rings. The
rings that we have already encountered have the property of being closed under
binomial, and not just polynomial, functions. We will need to formalise this
property. Note first that by \(2.2) there are unique polynomials $h^m_n(x)$,
$f_n(x,y)$ and $g^m_n(x)$, which are linear, bilinear and linear respectively
s.t.
\begin{eqnarray*}
{x \choose m}{x \choose n} &=&
h^m_n\left({x \choose 1},{x \choose 2},\dots ,{x \choose mn}\right)\\
{xy\choose n}&=&f_n\left( {x \choose 1},\dots ,{x \choose n},{y \choose 1},\dots ,
{y \choose n} \right)\\
{{x \choose m}\choose n} &=&g^m_n\left({x\choose 1},{x\choose 1},\dots ,
{x\choose mn} \right).
\end{eqnarray*}
\begin{definition}
A \Definition{numerical ring} is a commutative ring  $R$ together with 
functions  \pil{-\choose n}RR, $n\ge 0$, s.t.

\part[i]  ${0\choose n} = 1$,
                   
\part[ii] ${1\choose n}=0$, $n\ge 2$.

\part[iii] ${r\choose 1} = r$,

\part[iv] ${r+s\choose n}  =\sum _{i+j=n}{r\choose i}{s\choose i}$,

\part[v] ${r \choose m}{r \choose n} =
h^m_n\left({r \choose 1},{r \choose 2},\dots ,{r \choose mn}\right)$,

\part[vi] ${rs\choose n}=f_n\left( {r \choose 1},\dots ,{r \choose n},{s \choose
1},\dots , {s \choose n} \right)$,

\part[vii] ${{r \choose m}\choose n} =
g^m_n\left({r\choose 1},{r\choose 1},\dots ,{r\choose mn} \right)$.
\end{definition}
\begin{remark}
\part 
In the presence of \DHrefpart{v}, \DHRefpart{i}-\DHrefpart{iv} and
\DHRefpart{vi}-\DHrefpart{vii} are equivalent to $R$ being a $\lambda$-ring and
one can in fact replace the polynomials $f$ and $g$ by those used in the theory
of $\lambda$-rings (they are equal modulo \DHrefpart{vi}). Indeed, in the
presence of \DHrefpart{vi} the polynomials appearing in the theory of
$\lambda$-rings reduce to linear resp.~bilinear polynomials. As $f$ and $g$ are
characterised by \DHrefpart{iv} resp.~\DHrefpart{v} being true for $r,s\in
{\Z}$ and {\Z} is a special {\gla}-ring we see that they necessarily reduce to
$f$ and $g$.

\part In terms of the ring homomorphism \pil{\phi}R{1+R\bigl[[t]\bigr]} from the
theory of $\lambda$-rings the extra axiom \DHrefpart{v} can be described as
follows. The map $\phi$ can be thought of as giving an exponentiation of $1+t$
by elements of $R$ through $(1+t)^r:= \phi(r)$. For the exponentiation of an
arbitrary element $1+c(t)$ of $1+R\bigl[[t]\bigr]$ there are two
candidates. Either we can use the $R$-module structure on $1+R\bigl[[t]\bigr]$
given by $\phi$ or we can substitute $c(t)$ for $t$ in $(1+t)^r$. In the
presence of the $\lambda$-ring axioms, \DHrefpart{v} is equivalent to these two
constructions coinciding. The details are left to the reader.
\end{remark}
As usual we can construct for any set $S$ the free numerical ring ${\Num}(S)$ on
$S$, i.e.,\linebreak[4] $\Hom_{\Sets}(S,R)=\Hom_{\Num-rings}({\Num}(S),R)$ for any numerical ring
$R$ and also the free numerical ring ${\Num}^g(M)$ on an abelian group $M$. Then
${\Num}(S)={\Num}^g({\Z}[S])$.
\begin{lemma}\label{3.2}
For any set $S$, ${\Z}^S$ is a numerical ring with pointwise operations. Let
$x_j\in {\Num}_i$, $1\le j\le i$, be the projection on the $j$'th factor. Then
${\Num}_i = {\Num}(\{x_j: 1\le j\le i\})$.
\begin{proof}
It is clear that {\Z} with the binomial functions is a numerical ring, in fact
the axioms were set up precisely to ensure this. Hence ${\Z}^S$ certainly is a
numerical ring with pointwise operations. Furthermore, ${\Num}_i\subset
{\Z}^{\Z^i}$ is clearly stable under all operations and is hence a sub-numerical
ring. Therefore, if $S :=\{x_j\}$ there is a map of numerical rings ${\Num}(S)
\to {\Num}_i$ taking $x_j$ to $x_j$. By \(2.2) this map is surjective and if we
prove that
\begin{displaymath}
   {\Num}(S) = \sum  {\Z}{x_1\choose n_1}\dots {x_i\choose n_i}=:A,
\end{displaymath}
where the sum is not necessarily direct, then we are finished. As $x_j\in A$ it
is sufficient to show that $A$ is a numerical subring of ${\Num}(S)$. That $A$
is a subring follows from v) and stability under$-\choose n$ follows from the
rest of the axioms.
\end{proof}
\end{lemma}
For any complex $0\to C^0 \to C^1 \to C^2 \to \dots $ of abelian groups we may
construct a cosimplicial numerical ring as ${\Num}(C^\cdot
):={\Num}^g(\Gamma(C^\cdot))$, where ${\Num}(-)$ is extended pointwise to
simplicial objects. In case $C_i$ is f.g.~free for all $i$, then ${\Num}(C^\cdot
)=\Hom_{\Num}(\Gamma(\Hom_{\Z}(C^\cdot)),{\Z})$ giving the relation with the
preceding results.
\begin{lemma}\label{3.3}
Let  $R \to S$ and  $R \to T$ be morphisms of numerical 
rings. Then there is a structure of numerical ring on  $S\bigotimes_RT$ making 
it the  pushout, in the category of numerical rings, of  
$R \to S$ and  $R \to T$.
\begin{proof}
If $ S \to S\bigotimes_RT$ and $T \to S\bigotimes_RT$ are to be morphisms of
numerical rings then the definition of $-\choose n$ are forced by the axioms so
we begin by showing that they are well-defined. As was remarked above any
numerical ring is also a special $\lambda$-ring. This means that if $U$ is a
numerical ring and if we put
\begin{displaymath}
\funk{\phi\co U}{1+tU\bigl[[t]\bigr]}
        {r}        {\sum _{i=0}^\infty{r\choose i}t^i},
\end{displaymath}
then $\phi$ is a ring homomorphism where $1+tU\bigl[[t]\bigr]$ is given the ring
structure of \cite[Exp. V,2.3]{SGA6} with multiplication denoted by *. This
gives us ring homomorphisms $S \to 1+tS\bigotimes _RT\bigl[[t]\bigr]$ and $T \to
1+tS\bigotimes _RT\bigl[[t]\bigr]$ coinciding on $R$ and hence we get a ring
homomorphism $S\bigotimes _RT \to 1+tS\bigotimes RT\bigl[[t]\bigr]$ showing that
the operations $-\choose n$ are well-defined on $S\bigotimes _RT$. To show that
we get a numerical ring we reduce to $R={\Z}$ and $S$ and $T$ free numerical
rings on finite sets and conclude by \(2.2) and \(3.2) as these show that
$S\bigotimes T$ then is again a (free) numerical ring. Clearly $S\bigotimes _RT$
has the required universal property.
\end{proof}
\end{lemma}
We can now define twisted cartesian coproducts (TCcP's) of cosimplicial
numerical rings (CNR's) etc by dualising section 1 using the coproduct of lemma
\ref{3.3}.
\begin{definition}\label{fibred CNR}
A \Definition{fibred CNR} is a cosimplicial numerical ring  $Y$ 
s.\ t.\ $Y= \dli (\dots \to Y_n \to Y_{n+1} \to \dots )$ and

\part $\forall i\exists N:  Y^i_n \to  Y^i_{n+1}$ is an isomorphism for  $n\ge N$,

\part 
$Y_n \to Y_{n+1}$ is a PTCcP with cofibre of the form ${\Num}(\Gamma(C^\cdot _n
))$, where $C^\cdot _n$ is a bounded complex of free (not necessarily f.\ g.) 
abelian groups.
\end{definition}
We will say that a cosimplicial numerical ring $R$ is \Definition{connected} if
$H^0(R)={\Z}$ and $H^1(R)$ is torsion free and $1$-connected if $H^0(R)={\Z}$,
$H^1(R)=0$ and $H^2(R)$ is torsion free, where $H^*(R)$ denotes the cohomology
of the corresponding complex.
\begin{remark}
\part  In condition \DHrefpart{ii} we do not need the condition on the 
cohomology of $C^\cdot _n$; by refining and modifying the sequence 
$\dots \to Y_n \to Y_{n+1} \to \dots $ this condition can always be fulfilled.

\part The conditions defining connectivity and $1$-connectivity should be
considered in the light of the universal coefficient sequence; connectivity and
$1$-connectivity refers to vanishing of homology.
\end{remark}
\begin{theorem}\label{3.5}
\part[i] Let $X$ be a fibred CNR and $Y_1 \to Y_2$ a CNR-morphism which is a
cohomology equivalence. Then every morphism $X \to Y_2$ can be lifted to $Y_1$.

\part[ii]  Any cohomology equivalence between fibred CNR's is a 
homotopy equivalence. 

\part[iii] Let $Y$ be a $1$-connected CNR. Then there is a unique (up to homotopy)
fibred CNR $X$ and a numerical ring homomorphism $X \to Y$ which is a
cohomology equivalence.

\part Let  $Y$ be a connected CNR. Then there is a unique (up to 
homotopy) CNR  $X$ fulfilling the liftability with respect to cohomology
equivalences as in i) and a numerical ring 
homomorphism  $X \to Y$ which is a cohomology equivalence. 
\begin{proof}
\DHrefpart{i} is proved by successive liftings (and is essentially a numerical
cosimplicial version of the proof of the similar property for cdga's). The
lifting at one stage is accomplished as follows. By the definition of
\ref{fibred CNR} $C^\cdot_n$ will consist of normalised cochains in $X^n$ and
the subcomplex $K$ of $N(X^n)$ generated by $N(X^{n-1})$ and $C^\cdot_n$ is the
mapping cone of a map of complexes $N(X^{n-1})[1] \to C^\cdot_n$. We then extend
the lifting of $N(X^{n-1}) \to Y_2$ to $K$. This in turn gives a lift of
$\Gamma(C^\cdot_n) \to Y_2$ which then gives a lift of $\Num(\Gamma(C^\cdot_n))
\to Y_2$. Then \DHrefpart{ii} follows similarly. Let us turn to
\DHrefpart{iii}. We will build a Postnikov tower. Note to begin with that when
building this tower we must kill homology and not cohomology. Hence we assume
that we have \pil f{X^{n-1}}Y s.t.~if we consider $X^{n-1}$ and $Y$ as complexes
then $H^i(C(f))=0$ if $i<n$ and $H^n(C(f))$ torsion free. Let $C^0 \to C^1$ be a
free complex s.t.~$H^0(C^\cdot )=H^n(C(f))$ and $H^1(C^\cdot )={\rm
tor}H^{n+1}(C(f))$. Then there is a morphism of complexes $C^\cdot [-n] \to
C(f)$ inducing the identity on $H^n$ and the natural inclusion on
$H^{n-1}$. Hence there is a morphism of cosimplicial groups $\Gamma(C^\cdot
[-n]) \to \Gamma(C(f))$ and by composition a morphism $\Gamma(C^\cdot [-n-1])
\to X^n$. (Note that $C(f)$ fits in to a distinguished triangle $N(X^n) \to N(Y)
\to C(f)$.) The composite $\Gamma(C^\cdot [-n-1]) \to Y$ is, by construction,
nullhomotopic. By adjunction we get ${\Num}(C^\cdot [-n-1]) \to X^n$ whose
composite with $X^n \to Y$ is again nullhomotopic as adjunctions preserve
homotopies. Therefore there is a PTCcP $X^n \to X^{n+1} \to {\Num}(C\cdot [-n])$
and a lifting \pil g{X^{n+1}}Y. I claim that $H^i(C(g))=0$ for $i\le n$ and that
$H^{n+1}(C(g))$ is torsion free. To do this one has to say something about the
cohomology of ${\Num}(C^\cdot [-n])$. When $C^\cdot $ is finitely generated this
has already been done and the general case is done by approximating $C^\cdot $
by finitely generated subcomplexes. After that one has to look at the Serre
spectral sequence for the PTCcP constructed above. A small conceptual problem
arises as we want to kill homology but are working with cohomology. This can
certainly be overcome by brute force, but we will instead choose a hopefully
more conceptual approach. This entails, however, the introduction of
pro-(finitely generated abelian groups) and the reader who is unfamiliar with
the concept of pro-objects will have no problem in translating the proof to
follow into one using only cohomology. If $D^\cdot$ is a complex of torsion free
abelian groups then we define its homology by $H^i(D^\cdot) :=
"\ili"\{H_i(\Hom(D^\cdot _{\gal} ,{\Z}))\}$ (cf.~\cite[Exp. I,8]{SGA4:1}), where
$D^\cdot _{\gal}$ runs over all
finitely generated subcomplexes of $D^\cdot $. We have the usual universal
coefficient sequences expressing cohomology and homology in terms of each other
if we put, for an abelian group $M$, $\Hom(M,{\Z})$ resp.~$Ext^1(M,{\Z})$ equal to
$"\ili"\{\Hom(M_{\gal},{\Z})\}$ resp.~$"\ili"\{Ext1(M_{\gal},{\Z})\}$, where
$M_{\gal}$ runs over all f.g.~subgroups of $M$ and, for a pro-object
$\{M_{\gal}\}$, $\Hom(\{M_{\gal}\},{\Z})$ resp.~$Ext^1(\{M_{\gal}\},{\Z})$ equal
to $\dli \{\Hom(M_{\gal},{\Z})\}$ resp.~$\dli \{Ext^1(M_{\gal},{\Z})\}$. Hence
our assumptions imply that $H_i(C(f))=0$ if $i<n$ and we want to prove that
$H_i(C(g))=0$ if $i\le n$. Furthermore, we may present $C^\cdot $ as a direct
limit of complexes $C^\cdot _{\gal}$ which are finitely generated free,
concentrated in degrees $0$ and $1$ with $H^0$ free and $H^1$ torsion. Then
\begin{displaymath}
H_*({\Num}(C^\cdot [-n]))="\ili"\{H_*({\Num}(C^\cdot _{\gal} [-n]))\}
\end{displaymath}
 and by \(2.7:ii)        
\begin{displaymath}
 H_*({\Num}(C^\cdot _{\gal}[-n]))=H_*(K(H0(C^\cdot _{\gal}),n)).
\end{displaymath}
By the well-known computation of the cohomology of Eilenberg-MacLane spaces we
get that $\tilde H_i({\Num}(C^\cdot [-n]))=0$ if $i<n$ or $=n+1$ (as $n\ge 2$)
and $H_n({\Num}(C^\cdot [-n]))=H_0(C^\cdot )=H_n(C(f))$. Finally, as above we
get a Serre s.s. 
\begin{displaymath}
H_i(X^n,H_j({\Num}(C^\cdot [-n]))) \Rightarrow  H_{i+j}(X^{n+1}).
\end{displaymath}
This and the information we have on $H_j({\Num}(C^\cdot [-n]))$ gives an exact
sequence 
\begin{displaymath}
\diagram{
0&\to &H_{n+1}(X^{n+1})& \to &H_{n+1}(X^n)&\to &H_n(C(f))&\to \cr 
 &    &\nwarrow  & &\nearrow\cr
 &    &             &H_{n+1}(Y) \cr}
\end{displaymath}
\begin{displaymath}
\diagram{
\to &H_n(X^{n+1})& \to &H_n(X^n)&\to &0\cr
&\nwarrow&  &\nearrow\cr
&&H_n(Y)\cr}
\end{displaymath}
and isomorphisms $H_i(X^{n+1}) \to H_i(X^n)$ for $i<n$. By construction we have an 
exact sequence   
\begin{displaymath}
 H_{n+1}(Y) \to H_{n+1}(X^n) \to H_n(C(f)) \to H_n(Y) \to H_n(X^n)\to 0.
\end{displaymath}
Combining these two sequences we get that $H_{n+1}(Y) \to H_{n+1}(X^{n+1})$ is an 
epimorphism and that $H_i(Y) \to H_i(X^{n+1})$ is an isomorphism for $i\le n$ i.e.~that 
$H_i(C(g))=0$ for $i\le n$. In case $Y$ is only connected we only get that 
\begin{displaymath}
\coker(H_{n+1}(Y) \to H_{n+1}(X^{n+1})) \to  \coker(H_{n+1}(Y) \to H_{n+1}(X^n))
\end{displaymath} 
is zero so we may have to continue an infinite number of times just to kill
homology in one degree and the end result will not necessarily be a special
CNR. It will still have the lifting property of \DHrefpart{i} though.
\end{proof}
\end{theorem}
As should be no surprise there is a very tight relation between CNR's and
simplicial numerical objects.
\begin{proposition}
\part The functor that takes each CNR $X$ and associates to it the simplicial
scheme $\Sp X$ obtained by taking the spectrum in each degree is an equivalence
of categories between CNR's that are free finitely generated numerical ring in
each degree and a full subcategory $\Cal N$ of the category of simplicial
schemes. The inverse functor is taking global sections $\Gamma(X,\Cal O)$ of the structure sheaf.

\part The functor that to a simplicial scheme $Y$ associates its simplicial set
of $\Z$-points is an equivalence of categories from $\Cal N$ to the category of
simplicial numerical objects.

\part The functor that to a simplicial numerical scheme $X$ associates the CNR
$\Hom_{Rings}(X,\Z)$ induces an equivalence between the category of simplicial
numerical rings and the full subcategory of CNR's that are finitely generated
free in each degree.
\begin{proof}
The results clearly follow from \(2.2).
\end{proof}
\end{proposition}
We will use these equivalences to think of a numerical space as the $\Z$-points
$X(\Z)$ of a simplicial scheme $X$. Localisation and completion has a
particularly pleasant formulation in these terms; we have that $F(X(\Z) \Tensor
R)=X(R)$. In the case of completion we have the following rather striking fact
which in particular shows that $p$-complete homotopy types can be described in
terms of cosimplicial $\Z/p$-algebras whose cohomology is the cohomology of the
type. We also add a rather curious fact saying that in the $p$-complete case we
may use continuous chains to compute cohomology.
\begin{proposition}
\part Let $X \in \Cal N$. Then the reduction mod $p$ map $X({\Z}_p) \to X({\Z}/p)$ is
a bijection.

\part Assume $X$ is a $\Z_p$-numerical space and give the components of
$X(\Z_p)$ its $p$-adic topology. Then the cohomology of the complex of
$\Z_p$-continuous $\Z_p$-valued cochains is isomorphic with ordinary cohomology
of $X$ with $\Z_p$-coefficients.
\begin{proof}
The first part clearly amounts to showing that the reduction mod $p$ map\linebreak[4]
$\Hom_{Rings}({\Num}_i,{\Z}_p) \to \Hom_{Rings}({\Num}_i,{\Z}/p)$
is a bijection. This can no doubt be done directly but the ``real" reason why it
is true is the following. As we have seen, ${\Z}^i=Spec {\Num}_i$ is the Cartier
dual of $\mulf^i$, the product of $i$ copies of the formal multiplicative group,
and so \cite[II,\S4]{Dem72} for any ring $R$,
${\Z}^i(R)=\Hom_{\symb{\mbox{R-formal groups}}}(\mulf^i,\mulf)$ and it is well
known [loc.~cit.] that
$\Hom_{\symb{\mbox{R-formal grps}}}(\mulf^i,\mulf)={\Z}^i_p$ for $R={\Z}_p$
as well as ${\Z}/p$.

As for the second part we note that the ring of continuous functions from $\Z_p^n
\to \Z_p$ equals the $p$-adic completion of $\Num_n$, this is Mahler's theorem
(\cite{mahler58::aen}). Hence the complex of continuous cochains is the
completion of the complex of numerical cochains. This gives rise to short exact
sequences
\begin{displaymath}
\shex{\ili^1H^{i-1}_{\Num}(X,\Z/p^n)}{H^i_{cont}(X(\Z_p),\Z_p)}{\ili H^i_{\Num}(X,\Z/p^n)}
\end{displaymath}
but as the $H^i_{\Num}(X,\Z_p)$ are finitely generated $\Z_p$ the left hand side
is zero and the right hand side is $H^i_{\Num}(X,\Z_p)$.
\end{proof}
\end{proposition}
\begin{remark}
\part The first part of the proposition shows that the category of simplicial
$\Z_p$-numerical objects is equivalent to a category cosimplicial
$\Z/p$-algebras which in the case of $\Z_p$-numerical spaces computes the
$\Z/p$-cohomology of the space. This accords more with the usual view of
$p$-complete spaces where $\Z/p$-cohomology reflects isomorphisms. I do not
however know of an intrinsic characterisation of the algebras of the form
$\Num_i/p$ in the style of characterising $\Num_i$ as free numerical
algebras. It should be noted that there is a Stone type duality between
$\Z/p$-algebras $R$ fulfilling $r^p=r$ for each $r \in R$ and totally
disconnected compact spaces; the space is the set of ring homomorphisms into
$\Z/p$ and the ring is the set of continuous maps into $\Z/p$. For $p=2$ this is
the usual Stone duality.

\part By \cite[V,Thm 2.3.10]{lazard65::group} we get that the cohomology of a
$p$-complete finitely generated torsion free nilpotent group can also be
computed using analytical cochains.
\end{remark}
Let us end this section with an observation that shows that disregarding the
rest of the section cosimplicial numerical rings are related to homotopy
theory. Thus let $R$ be a cosimplicial numerical ring and $\alpha \in
H^i(R)$. If $z$ is a representing cocycle in $N(R)$ for $\alpha$ it is
represented by a map $\Z[-i] \to N(R)$ and hence a map $\Gamma(\Z[-i]) \to R$
and again by a map $\Num(\Gamma(\Z[-1])) \to R$. This induces a map on
cohomology $H^*(\Num(\Gamma(\Z[-1]))) \to H^*(R)$. By \(2.7) this means that all
cohomology operations will operate on the cohomology of cosimplicial numerical
rings with all relations being preserved.
\end{section}

\begin{section}{Nilpotent groups}

We will spend some time considering the case of $K(G,1)$'s or equivalently
nilpotent groups. More precisely we will only consider those that are torsion
free. It can be concluded from the results of the previous section that each such
group $G$ may be identified as a set with $\Z^n$ for some $n$ in such a way that
the multiplication and inverse are given by numerical maps and that the
cohomology may be computed using numerical cochains. We will now see that there
is a \emph{canonical} way to define the structure of object in $\Num$ on the set
underlying a nilpotent group such that the group structure is given by a group
object in $\Num$. Indeed, let $G$ be a nilpotent f.\ g.\ torsion free group and let
$\Z[G]$ be its group algebra. Any function \pil{\phi}G{\Z} gives rise to an
additive function, also denoted $\phi$ \pil{\phi}{\Z[G]}{\Z} using the fact that
$\Z[G]$ is free on $G$. We say that $\phi$ is \Definition{$\cdot$-numerical}, where
$\cdot$ is the product on $G$, if $\phi$ vanishes on some power of the
augmentation ideal of $\Z[G]$. We denote by $\Num_{\cdot}(-,\Z)$ the set of $\cdot
$-numerical functions. 
\begin{remark}
These functions were first introduced by Passi (\cite{passi68::dimen}) and are
also known as ``Passi polynomial'' maps. We have chosen a different terminology
because $+$-numerical functions on $\Z^n$ are exactly numerical functions and
because of subsequent results.
\end{remark}
We need a preliminary result giving a characterisation of
$\cdot$-numerical functions that may be of independent interest. For that recall
that a module for a group $G$ is said to be \Definition{unipotent} if it is a
successive extension of modules with trivial action.
\begin{lemma}\label{numAndUni}
Let $(G,\cdot)$ be a torsion-free f.\ g.\ nilpotent group. Then a map $G \to \Z$
is $\cdot$-numerical if and only if it generates a $\Z$-finitely generated
unipotent submodule of $\Z^G$.
\begin{proof}
It is clear that a $G$-module $M$ is unipotent if and only if it is annihilated
by some power of the augmentation ideal of the group ring $\Z[G]$. By definition
a map $G \to \Z$ is $\cdot$-numerical if and only if it is annihilated by a
power of the augmentation ideal and as the augmentation ideal is two-sided this
is true if and only if it generates a submodule that is. Finally, as $\Z[G]$
modulo any power of the augmentation ideal is a finitely generated $\Z$-module,
any unipotent submodule of $\Z^G$ generated by one element is finitely generated.
\end{proof}
\end{lemma}
We will need the following result seemingly unrelated result. Recall that a
(smooth) connected algebraic group is \Definition{unipotent} if it is an
algebraic subgroup of the group of unipotent upper triangular $n\times
n$-matrices for some $n$ and that in that case every linear representation of
it is unipotent (cf.\ \cite[Cor.\ 3.4]{raynaud70::group}). In particular the
points of the group over a base field is a nilpotent group. Furthermore, if that
base field is the rational numbers, if $g$ is a point defined over the it and
$f$ is a polynomial vanishing on the group then we may introduce the polynomial
in $x$ $g^x := \exp(x\log(g))$, where the logarithm is a finite series mapping
unipotent upper triangular matrices to nilpotent ones and the exponential is a
finite series mapping nilpotent upper triangular matrices to unipotent
ones. This polynomial vanishes on all integers and hence is identically zero. In
particular it vanishes on $g^r$ for $r \in \Q$. As the group is closed in the
Zariski topology it is equal to the common zero set of all such $f$ and $g^r$ is
in the group. This means that the group of rational points is a uniquely
divisible nilpotent group.
One of the most natural questions on the relation between nilpotent torsion-free groups
and numerical groups is answered by the following result.
\begin{lemma}\label{numIsNil}
Let $G$ be a group object in $\Num$. Then the underlying group is a finitely
generated torsion-free nilpotent.
\begin{proof}
If we extend the scalars of $\Hom_{\Num}(G,\Z)$ to $\Q$ we get a polynomial ring
over $\Q$ which is the affine algebra of an algebraic group $\sG$ over
$\Q$. Hence \cite{lazard55::nilp} applies and we conclude that $\sG$ is unipotent and
hence that $\sG(\Q)$ is a uniquely divisible nilpotent group. Clearly, $G$ is a
subgroup of $\sG(\Q)$ and hence is nilpotent and torsion-free. To prove finite
generation we note that there as a finite dimensional faithful subrepresentation
$V$ of the representation of $\sG$ on its affine algebra. We now want show that
there is finitely generated subgroup $M$ of $V$ stable under the action of
$G$. As $V$ is finite dimensional it contains a finite number of vectors
spanning it as $\Q$-vector space. After possibly multiplying them by a non-zero
integer we may assume that they are contained in $\Hom_{\Num}(G,\Z)$. It is
therefore sufficient to show that each $f \in \Hom_{\Num}(G,\Z)$ lies in a
finitely generated subgroup of $\Hom_{\Num}(G,\Z)$ invariant under $G$. For this
we consider the product map \pil\varphi{G\times G}G and write the pullback
$\varphi^*f \in \Hom_{\Num}(G,\Z) \Tensor \Hom_{\Num}(G,\Z)$ as $\sum_if_i \tens
g_i$. This means that for $g,h \in G$ $f(g\cdot h)=\sum_i f_i(g)g_i(h)$ and by
definition $h \in G$ acts on $f$ by $(hf)(g)=f(gh)$. Hence we get that $hf =
\sum_i g_i(h)f_i$ so that the translates of $f$ by the elements of $G$ lies in
the finitely generated group spanned by the $f_i$ and hence is finitely
generated.

Thus, $G$ is a subgroup of the subgroup $\sG_M$ of elements of $\sG(\Q)$ stabilising
$M$. As a subgroup of a finitely generated nilpotent group is finitely generated
it is enough to show that $\sG_M$ is finitely generated. This will be done by
induction over the dimension of $V$ (with $\sG$ changing during the
induction). As $\sG$ is unipotent $V$ contains a $1$-dimensional subspace $U$ on
which $\sG$ acts trivially. We may use the induction hypothesis on the image of
$\sG$ in $\Aut(V/U)$ and the image $M'$ of $M$ in $V/U$ to conclude that the image of
$\sG_M$ in $\Aut(V/U)$ and it is then enough to show that the kernel of this map
is finitely generated. However, that kernel is a subgroup of the abelian group
of additive maps $\Hom(M',U\cap M)$ which is finitely generated.
\end{proof}
\end{lemma}
\begin{proposition}
Let $(G,\cdot)$ be a finitely generated torsion-free nilpotent group.

\part[0] $G$ may be identified with $\Z^n$ in such a way that multiplication and 
inverse on $G$ are numerical functions and $K(G,1)$ is a numerical space.

\part[i] $\Num{\cdot}(G,\Z)$ is a numerical subring of the numerical ring of all
functions $G\to \Z$.

\part[ii] If $G$ has been given the structure of group object in $\Num$ for
which $K(G,1)$ is a numerical space then a function $G \to \Z$ is
$\cdot$-numerical precisely when it is numerical with respect to the given
numerical structure.

\part[iii] The numerical ring $\Num{\cdot}(G,\Z)$ is isomorphic to the free numerical
ring on a finite number of generators and the natural map $G\to
\Hom(\Num{\cdot}(G,\Z),\Z)$ is a bijection.

\part[iv] The product map $G\times G\to G$ induces a map $\Num{\cdot}(G,\Z)\to
\Num{\cdot}(G,\Z)\bigotimes \Num{\cdot}(G,\Z)$, where $\Num{\cdot}(G,\Z)\bigotimes
\Num{\cdot}(G,\Z)$ is thought of as a subring of the set of functions $G\times G\to
\Z$. The inverse $G\to G$ induces a map $\Num{\cdot}(G,\Z)\to\Num{\cdot}(G,\Z)$.
\begin{proof}
\DHrefpart{0} follows from lemma \ref{2.7} and induction over the length of the
ascending central series,. \DHrefpart{i} is obvious. As for \DHrefpart{ii}
consider first a function \pil fG{\Z} that is numerical with respect to the
given numerical structure. Let $\sG$ be the algebraic group over $\Q$ whose ring
of regular function is $\Hom_{\Num}(G,\Z)\Tensor\Q$. Then $f$ generates a finite
dimensional subrepresentation of $\Hom_{\Num}(G,\Z)\Tensor\Q$ which is unipotent
as $\sG$ is by \cite{lazard55::nilp}. Hence $f$ $\cdot$-numerical by
\(numAndUni). Assume conversely that \pil fG{\Z} is $\cdot$-numerical. Again by
lemma \ref{numAndUni} it generates a unipotent module $M$. We may choose a
$G$-invariant filtration of $M$ whose successive quotients are free of rank $1$
with trivial $G$-action. Having done this, the $G$-action on $M$ corresponds to
a group homomorphism from $G$ to $U$, the group of unipotent upper triangular
integer $n\times n$-matrices, where $n$ is the rank of $M$. Furthermore, $U$ is
a numerical group (in fact an algebraic one) and $f$ is the composite of a
numerical map $U \to \Z$ and the group homomorphism $G \to U$. It will therefore
suffice to show that the group homomorphism $G \to U$ is numerical. The
ascending central series of $U$ is given by $\{U_m\}$, where $U_m$ is defined by
\set{(a_{ij})}{a_{ij}=0\; \mathrm{if}\; j < i <= j+m} and $U/U_m$ is clearly
also a numerical group. We now prove by descending induction on $m$ that the
composite $G \to U \to U/U_m$ is numerical. The homomorphism $U/U_m \to
U/U_{m-1}$ is a central extension. The obstruction for lifting the homomorphism
$G \to U/U_{m-1}$ to $U/U_m$ is an element of $H^2(G,U_{m-1}/U_m)$ that is zero
as the morphism is known to lift and the obstruction for lifting it to a numeric
homomorphism is an element of $H^2_{\Num}(G,U_{m-1}/U_m)$. The latter group maps
bijectively to the former by \(2.7) and hence the numerical homomorphism $G \to
U/U_{m-1}$ lifts to a numerical homomorphism $G \to U/U_m$. The set of group
homomorphism liftings are classified by $H^1(G,U_{m-1}/U_m)$ and the set of
numerical group homomorphism liftings are classified by
$H^1_{\Num}(G,U_{m-1}/U_m)$. Again by \(2.7) the map between these groups is a
bijection and hence every lifting is numerical. In particular the given one is
which finishes the induction step.

Finally, \DHrefpart{iii} and \DHrefpart{iv} follow from \DHrefpart{0} and \DHrefpart{ii}.
\end{proof}
\end{proposition}
We gather together the main results of this section in the following theorem.
\begin{theorem}
\part The abstract group underlying a group object in $\Num$ is a finitely
generated torsion free nilpotent group.

\part The forgetful functor from the category of group objects of $\Num$ to the
category of finitely generated torsion free nilpotent groups is an equivalence
of categories.

\part[iii] Let $G$ be a finitely generated torsion-free nilpotent group and let
$S$ be the ring of $\cdot$-numerical functions on $G$. Then $S$ is a free
numerical ring on the rank of $G$ generators. The product, inverse and unit
element of $G$ induces a Hopf algebra structure on $S$, the evaluation map $G
\to \Hom_{Rings}(S,\Z)$ is a bijection and through this bijection, the Hopf
algebra structure on $S$ induces the given group structure on $G$.

\part With notations as in the previous part, the group cohomology $H^*(G,\Z)$
of $G$ may be computed using $\cdot$-numerical cochains.

\part With notations as in \DHrefpart{iii} let $R$ be a coefficient
ring. Then the group structure induced on $G_R := \Hom_{Rings}(S,R)$ by the Hopf
algebra structure of $S$ together with the group homomorphism $G \to G_R$ given
by the composite of the isomorphism $G \to \Hom_{Rings}(S,\Z)$ given by
\DHrefpart{iii} and the map induced by the inclusion $\Z \to R$ is an
$R$-localisation.
\begin{proof}
This follows from the previous results of this section together with \(localisation).
\end{proof}
\end{theorem}
Apart from possible group theoretic applications we can apply our results to
more general homotopy types.
\begin{proposition}
Let $G_\cdot$ be a simplicial group all of whose components $G_n$ are
f.g.~torsion free nilpotent groups. Then $G_\cdot $ has a natural structure of
simplicial group in $\Num$ and using this structure $K(G_\cdot )$, the simplicial
classifying space of $G_\cdot $, has a natural structure of simplicial object in
$\Num$. As such it is a numerical space.
\begin{proof}
This follows directly from  the spectral sequence
\begin{displaymath}
E_1^n = H^*_\Z(K(G_n,1)) \Rightarrow  H^*_\Z(K(G_\cdot ,1))
\end{displaymath}
and the corollary.
\end{proof}
\end{proposition}
This result gives a relation between the present approach and one given by
Quillen to rational homotopy theory (cf.~\cite{Qu69}). Indeed, there
Quillen represents any finitely generated complex up to homotopy by exactly a
$K(G_\cdot )$ as in the proposition. What the proposition shows is that it also
gives a representation of the complex as a numerical space. The next step
in Quillen's construction is to pass to the Malcev completion of the components
of $G_\cdot$ which fits very well in our context as taking the $\Q$-points of
the $G_n$'s considered as group schemes.
\end{section}
\begin{section}{Sullivan models}

We now want to see how Sullivan's theory (cf., \cite{sullivan77::infin}) of
minimal models fits in with the present theory. His theory is a rational so
throughout this section the coefficient ring will be the ring of rational
numbers. As numerical functions then are the same as polynomial ones we will
call them just that.  A numerical space then can be seen as a simplicial scheme
which in each degree is an affine space. Let us introduce some notation
appropriate to the situation. We let $\Delta_a$ be the cosimplicial scheme that
in each degree $n$ is the algebraic $n$-simplex $\Sp
\Q[x_0,\dots,x_n]/(x_0+\dots+x_n-1)$ with the obvious face and degeneracy
operators and let $\Omega_a$ be the simplicial graded commutative differential
graded algebra (\Definition{cdga}) of algebraic forms on $\Delta_a$. Let us
recall that Sullivan associates to each simplicial set $X$ the differential
graded algebra $\Cal E(X)$ consisting of a choice of forms on
$X_n\times\Delta_a^n$ with appropriate compatibility conditions with respect to
face and degeneracy operations, where $X_n$ is thought of as a zero-dimensional
scheme being the disjoint union of copies of $\Sp\Q$, one for each point of
$X_n$. To generalise this to the case of a numerical space $X$ we consider
relative forms on $X_n\times\Delta_a^n$, relative to the projection on the first
factor, thus obtaining a cdga $\Cal E_a(X)$. Another way to think of this is to
consider the set of simplicial algebraic maps from $X$ to $\Omega_a$, where an
algebraic map from an affine space $Y$ over $\Q$ to a $\Q$-vector space $V$ is a
map $Y \to V$ whose image lies in a finite dimensional subspace $U$ of $V$ and
is algebraic as a map $Y \to U$. Now, one proves as in
\cite[Thm.~7.1]{sullivan77::infin} that the cohomology of $\Cal E_a(X)$ computes
the numerical cohomology of $X$. In particular, if $X$ is a numerical space the
inclusion map $\Cal E_a(X) \to \Cal E(X)$ is a quasi-isomorphism. Note that even
though $\Cal E_a(X)$ is considerably smaller than $\Cal E(X)$ even when $X$ is a
minimal numerical space $\Cal E_a(X)$ is far from being a minimal model (or even
a model), in fact in general $\Cal E_a(X)$ will not be connected (i.\ e., $\Q$
in degree $0$). It would interesting to have a modification of this construction
that would give a minimal model from a minimal numerical space\dots
\begin{remark}
Instead of looking at relative forms on $X_n\times\Delta_a^n$ one could look at
all forms. This would give a complex analogous to the de Rham complex of a
simplicial manifold and its cohomology does in fact compute the algebraic de
Rham cohomology of the simplicial scheme $X$. However, as each component $X_n$
is contractible this de Rham complex is acyclic.
\end{remark}
The relation with Sullivan's geometric realisation functor (cf.,
\cite[\S8]{sullivan77::infin}) seems to be more interesting. By investigating it
a little bit more closely than is done (\cite{sullivan77::infin}) we will find a way
of giving a direct construction of a minimal numerical space from a minimal
cdga. For that we need some preliminary results. We begin by recalling that the
\Definition{canonical truncation}, $\tau_{\le n}$, of a complex $C^\cdot$ is the
subquotient of $C^\cdot$ obtained by first taking the quotient by the subcomplex
$C^{>n}$ and then the subcomplex which is unchanged in degrees $<n$ and the kernel
of the differential \pil{d^n}{C^n}{C^{n+1}} in degree $n$. Normally $\tau_{\le
n}C^\cdot$ is seen as a subcomplex of $C^\cdot$ but done in this fashion it is clear
that if $C^\cdot$ is a cdga then we get an induced structure of cdga on $\tau_{\le
n}$. Note that the cohomology groups of $\tau_{\le n}C^\cdot$ are the same as those
for $C^\cdot$ in degrees $\le n$ and $0$ otherwise.
\begin{definition-lemma}
\part For each $n \ge 0$ the simplicial group $\Omega^n$ is acyclic.

\part Putting $Z^n := \ker d: \Omega^n \to \Omega^{n+1}$, the simplicial abelian
group $Z^n$ is $0$ in degrees $< n$, $\pi_n(Z^n)=\Q$ and $\pi_i(Z^n)=0$ if $i >
n$. As $Z^n_i=0$ for $i > n$ there is a mapping \pil{\int}{Z^n}{K(\Q,n)}
inducing the identity on $\pi_n$. We let $T^n$ be the simplicial cdga that is
the quotient of $\tau_{\le n}\Omega$ by the kernel of $\int$.
\begin{proof}
The first part is implicit in \cite{sullivan77::infin} but can be found
explicitly in \cite[Lemma 10.7 \& \S17, Ex.~3]{felix01::ration}. The second part
is clear when $n=0$ as $Z^0$ is the constant simplicial object with constant
value $\Q$. For larger $n$ it follows from the first part, induction and the
exact sequences
\begin{displaymath}
\shex{Z^n}{\Omega^n}{Z^{n+1}}
\end{displaymath}
coming from the acyclicity of $(\Omega,d)$.
\end{proof}
\end{definition-lemma}
\begin{remark}
It is easily seen that the identification $\pi_n(Z^n)$ with $\Q$ can also be
given by the map $Z^n_n \to \Q$ given by integration over the standard simplex
hence justifying the terminology.
\end{remark}
Recall now, (\cite[\S8]{sullivan77::infin}), that the spatial realisation of a
commutative differential graded $\Q$-algebra $A$, that we will assume is locally
finite, i.\ e., finite-dimensional in each degree, is the simplicial set $\spres
A$ that in degree $n$ is the set of cdga-maps from $A$ to $\Omega_n$. The set of
such maps is in a natural fashion the $\Q$-points of a formal $\Q$-scheme. More
precisely, it is a closed subset of the formal affine space of graded linear
maps from $A^{\le n}$ to $\Omega_n$ (formal as the target space is infinite
dimensional when $n >0$). Hence, the spatial realisation is a simplicial formal
scheme. As will be seen in a moment it is very big when $A$ is a (minimal) model
but we will want to cut it down to reasonable size. For two $m$-simplices
\pil{f,g}A{\Omega} we define equivalence relations for each $n > 0$
\begin{condition}{$f \sim_n g$}
$\iff$\hfil\parbox{0.9\textwidth}{The two induced maps \pil{f',g'}{\tau_{\le
n}\gamma_{n}A}{\tau_{\le n+1}\Omega}, where $\gamma_{\le n}A$ denotes the sub-cdga of $A$
generated by the elements of degree $\le n$, coincide when composed by the
surjection $\tau_{\le n+1}\Omega \to T^{n+1}$.}
\end{condition}
We then define a quotient $\sprest A$ of $\spres A$ defined by the intersection
of all these equivalence relations. Clearly $\sprest A$ is contravariantly
functorial in $A$ and the quotient map $\spres - \to \sprest -$ is a natural
transformation.
\begin{proposition}\label{t-realisation}
Let $A$ be a locally finite cdga, $V$ a graded finite dimensional vector space
concentrated in degree $n > 0$ and $A\Tensor\Lambda V$ a cdga which as a graded
algebra is the graded tensor product of $A$ and the free graded commutative
algebra on $V$ and such that $d$ maps $V$ into $A \Tensor\Q$. By functoriality
the inclusion mapping $A \to A\Tensor\Lambda V$ then induces maps
$\spres{A\Tensor\Lambda V} \to \spres A$ and $\sprest{A\Tensor\Lambda V} \to
\sprest A$.

\part $\spres{A\Tensor\Lambda V} \to \spres A$ is a PTCP of
simplicial formal schemes with structure group $\Hom(V,Z^n)$.

\part  $\sprest{A\Tensor\Lambda V} \to \sprest A$ is a PTCP of
simplicial formal schemes with structure group\linebreak[4] $K(\Hom(V,\Q),n)$. The quotient
mapping $\spres{A\Tensor\Lambda V} \to \sprest{A\Tensor\Lambda V}$ is a mapping
of PTCP's over $\spres A \to \sprest A$ with respect to the structure group map
induced by the natural surjection $\Z^n \to K(\Q,n)$.
\begin{proof}
If we begin with the first part and we consider a fibre of
$\spres{A\Tensor\Lambda V} \to \spres A$ over an $m$-simplex
\pil{\phi}{A}{\Omega_n}, an extension to $A\Tensor\Lambda V$ is completely
determined by the restriction of $d$ to $V$. Furthermore, $d$ applied to any
element $v$ of $V$ is determined as it has to be the already prescribed image of
$dv \in A$. That means that if $f$ and $g$ are two extensions of $\phi$ then
$f-g$ maps $V$ into $Z^n_m$. Conversely given an extension $f$ and a linear map
$V \to Z^n_m$ there is an extension $g$ of $\phi$ such that the restriction of
$f-g$ to $V$ is the given map. Hence, the map $\spres{A\Tensor\Lambda V} \to
\spres A$ is a principal homogeneous space over the simplicial group $Z^n$. To
show that it is a PTCP we need to find a section the restriction
$\spres{A\Tensor\Lambda V} \to \spres A$ to the subcategory $\Delta_*$ of
$\Delta$ consisting of those increasing maps $\{0,1,\dots,m\} \to
\{0,1,\dots,n\}$ that take $0$ to $0$ (cf., \cite[18.7]{May92}). This is
obtained by the following observations
\begin{itemize}
\item Given a map \pil{\phi}A{\Omega_m} to extend it to $A\Tensor\Lambda V$ one
needs to find a map \pil f{V}{\Omega_m} such that $df(v)=\phi(dv)$ for all $v
\in V$. This is possible as $d\phi(dv)=0$ and $\Omega_m$ is acyclic. It can be
done explicitly given a contraction of $\Omega$.

\item Given a $\Q$-point $b \in \Delta_a^m$ we may use it as origin and use the
algebraic contraction $x \to t(x-b)+b$ and the usual integration formulas to
construct a contraction of $\Omega_m$. This contraction is natural for affine
maps preserving the chosen basepoints.

\item The restriction of the cosimplical scheme $\Delta_a$ to $\Delta_*$ has a
base point and hence the restriction of $\Omega$ to $\Delta_*$ has a contraction.
\end{itemize}

To turn to the second part it is clear that the only of the equivalence
relations $\sim_m$ on $\spres{A\Tensor\Lambda V}$ that does not factor through
$\spres{A\Tensor\Lambda V} \to \spres A$ is $\sim_n$. As for $\sim_n$ it is
clear that two  $m$-simplices $f$ and $g$ are equivalent if their restrictions to $A$
are and if the restrictions to $V$ of $\int\circ d\circ f$ and $\int\circ d\circ
g$ are equal. This combined with the first part now gives the second.
\end{proof}
\end{proposition}
From the proposition the main result of this section immediately follows.
\begin{theorem}
Let $A$ be a nilpotent cdga model. Then the natural map $\spres
A\to\sprest A$ is a homotopy equivalence and $\sprest A$ has a natural structure
of special $\Q$-numerical space. It is minimal if $A$ is.
\begin{proof}
We leave to the reader to prove, in a fashion analogous to
\cite{sullivan77::infin}, that if $A' \to A$ is a minimal model then $\sprest A
\to \sprest{A'}$ is a homotopy equivalence and that the formal scheme structure
on $\spres A$ induces a special numerical space structure on $\sprest A$ and
will assume that $A$ is minimal. That means that there is a filtration $A^n$ of
$A$ by sub-cdga's such that $A^n = A^{n-1}\Tensor\Lambda V_n$, where $V_n$ is
concentrated in a single degree and $d$ maps $V_n$ into $A^{n-1}$. Furthermore,
the degree of $V_n$ tends monotonically to infinity with $n$. One now proves by
induction that $\spres{A^n} \to \sprest{A^n}$ is a homotopy equivalence and that
$\sprest{A^n}$ is a minimal $\Q$-numerical space using proposition
\ref{t-realisation}. One then concludes by noting that $\sprest A$ is the
inverse limit of $\dots\to\sprest{A^n} \to \sprest{A^{n-1}}\to\dots$ and that
this system is eventually constant in each degree.
\end{proof}
\end{theorem}
\end{section}

\begin{section}{Fibrations}

As a simple, and not very original, application of the ideas of this paper we
will study fibrations.  We now note that we can relativise our constructions;
a morphism $R \to S$ of cosimplicial numerical rings may be factored $R \to S'
\to S$ where $R \to S'$ is a direct limit of a succession cTCP's with fibers as
before and $S' \to S$ is a cohomology equivalence. We call such a factorisation
a \Definition{special resolution} of the map. We get as before that any two
special resolutions are homotopic and for any map $R \to T$ of cosimplicial
numerical rings we will call $T \to T\bigotimes _RS'$ the homotopy pushout of $R
\to S$.
\begin{definition}
Let \pil{\phi}RS be a morphism of cosimplicial numerical rings. We say that
$\phi$ is a cofibration if for one (and hence any) special resolution $R \to S'
\mapright\rho S$ of $\phi$ and every morphism $R \to T$ of cosimplicial
numerical rings, $\rho\otimes _RT$ is a cohomology equivalence.
\end{definition}
We then have the following result.
\begin{proposition}
\part Any morphism  $R \to S$ of cosimplicial numerical 
rings which is flat in each degree (i.e.~$R_n \to S_n$ is a flat 
ring homomorphism for each  $n$) is a cofibration.

\part Let \pil fXY be a morphism of numerical spaces. The 
homotopy pullback of the map of simplicial sets underlying  $f$ by 
any map of numerical spaces has cohomology equal to the homotopy 
pushout of the corresponding cosimplicial numerical rings of 
numerical functions.

\part If  $R \to S$ is a fibration and  $R \to T$ a morphism, then there is a 
spectral sequence
               
\begin{displaymath}
Tor_*^{H^*(R)}(H^*(S), H^*(T)) \Rightarrow   H^*(S\bigotimes _RT).
\end{displaymath}
\begin{proof}
 i) is clear as a special morphism is flat so $S\bigotimes _R(-)$ and 
$S'\bigotimes _R(-)$ are both exact. As for ii) one verifies it by induction over a 
Postnikov tower of $X \to Y$. Finally, iii) follows as in the simplicial case 
\cite[II, Thm 5]{Qu67}.
\end{proof}
\end{proposition}
\begin{remark}
In the case of a diagram of spaces this gives the Eilenberg-Moore spectral 
sequence.
\end{remark}
\end{section}

\bibliography{abbrevs,algtop,alggeom,arithm,algebra}
\bibliographystyle{pretex}

\end{document}